\documentstyle[12pt]{article}
\newcommand{\p}{\it{p}}
\newcommand{\h}{\it{h}}
\newcommand{\g}{\it{g}}
\newcommand{\ga}{\it{a}}

\newcommand{\k}{\it{k}}
\newcommand{\n}{\it{n}}
\newcommand{\m}{\it{m}}
\newcommand{\mibitem}[1]{\bibitem{#1}}

\begin{document}
\begin{titlepage}
\begin{flushright}
UU-ITP 09-1995 \\
\end{flushright}

\vskip 0.3truecm

\begin{center}
{ \large \bf

REMARKS ON GEOMETRIC QUANTIZATION
\\ \vskip 0.2cm
OF R-MATRIX TYPE POISSON BRACKETS}
\\
\end{center}

\vskip 3.0 cm

\begin{center}

{ \large \bf  Alexey Kotov $^{*\dagger}$ } \\

\vskip 0.6cm

{\it Department of Theoretical Physics, Uppsala University \\

P.O. Box 803, S-75108, Uppsala, Sweden  }\\

\vskip 0.4cm

\end{center}

\vskip 4.0cm

\rm
\noindent
  We check the Vaisman condition of geometric quantization for
R-matrix type Poisson pencil on a coadjoint orbit of a compact
semisimple Lie group. It is shown  that this condition isn't
satisfied
for hermitian symmetric spaces . We construct also some examples when
Vaisman condition
takes place.

\vfill

\vskip 0.4cm

\begin{flushleft}
\noindent
\rule{5.1 in}{.007 in}\\
$^{*}${\small E-mail: $~~$ {\small \bf kotov@rhea.teorfys.uu.se}  \\
$^{\ddagger}$  permanent address : Chair of Nonlinear Dynamic Systems and
Control Processes, Department of Computational Mathematics and
Cybernetics, Moscow State University, Vorob'evy Gory, Moscow 119899, Russia

IMEM, Odessa State University, Petra Velikogo 2, Odessa 270000, Ukraine}\\

\end{flushleft}

\end{titlepage}

\vfill\eject

\baselineskip 0.65cm

\section{Introduction}

 Let  $G$ be a semisimple Lie group and  $\it r$ be  the {\em
Drinfeld-Jimbo R-matrix} which satisfies the {\em modified
Yang-Baxter}  equation.

Consider an orbit  $\cal O$ in  $\g^{\ast}$          of the coadjoint
action of $G$.
  There exists a Poisson bracket  $\{\;,\;\}_{KKS}$  (the {\em
Kirillov-Kostant-Souriau} bracket) on $\cal O.$

   We can introduce another Poisson structure on $\cal O$ by two
ways.\\
  (i)  $G$ acts on $\cal O$ and preserves the KKS-bracket, therefore
there is a representation
\smallskip
   \begin{equation}
   \rho : {\g } \rightarrow Der \left(C^{\infty}
(M),\{\;,\;\}_{KKS}\right)
   \end{equation}

\smallskip
 The map  $\{\;,\;\}_{\it r } :  C^{\infty} (M)^{\otimes 2} \to
C^{\infty} (M)$
   given by $\{f_1,f_2\}_{\it r }= \langle (\rho \otimes \rho ) {\it
r} , (df_1\otimes df_2) \rangle$ defines a skewsymmetric
multiplication on $C^{\infty} (M)$. It satisfies the Jacobi identity
iff orbit  $\cal O$ is a R-matrix type. The restriction of  the
3-vector  $[ [{\it r,r}] ]= [ {\it r}_{12},{\it r}_{13}]+[{\it
r}_{12},{\it r}_{23}]+[{\it r}_{13},{\it r}_{23} ]$, where ${\it
r}_{12}={\it r} \otimes1$, on $\cal O$ is identically zero in this
case. It's true if and only if  $\cal O$ is an orbit of  a semisimple
element and $\cal O$ is a symmetric space or  $\cal O$ is an orbit of
nilpotent element of height 2$\bf [3]$.The brackets $\{\;,\;\}_{\it r}$ and
$\{\;,\;\}_{KKS}$ are always compatible or form a {\em Poisson
pencil} on this orbit.

\smallskip

  (ii) Let  $K$  be a compact real form of  the semisimple complex
Lie group  $G$  and  $\cal O$  is an orbit in  ${\k}^{\ast}$.
   Let us consider  $\cal O$ as a Poisson coset space of Poisson Lie
group  $G$  with the coboundary {\em Sklyanin-Drinfeld} brackets
defined on  $G$
 as follows
\smallskip
\begin{equation}
\{f_1,f_2\} = \{f_1,f_2\}_L - \{f_1,f_2\}_R = \langle
(\rho_L^{\otimes 2} - \rho_R^{\otimes 2} ) ( {\it r} ), df_1 \otimes
df_2 ) ,
\end{equation}
\smallskip
 where  $\rho_R ( \rho_L )$     is the representation of the
corresponding Lie algebra  $\g$  in the space  $C^{\infty} (M) $ by
the left-(right-) invariant vector fields .\\
  The proper structure on  $\cal O$ can be obtained as a result of
Poisson reduction  $\bf [9]$. It is also called the Sklyanin-Drinfeld
brackets and denoted as  $\{\;,\;\}_{SD}$.
  The Kirillov-Kostsant-Souriau and Sklyanin-Drinfeld brackets form a
Posson pencil iff $\cal O$ is a hermitian symmetric space  $\bf
[6]$.\\
  It can be proved that the Poisson pairs   $a\{\;,\;\}_{KKS} +
b\{\;,\;\}_{SD}$  and  $a\{\;,\;\}_{KKS}+b\{\;,\;\}_{\it r}$
coincide on hss.

  The simple inspection of this Poisson  structures associated with
$\it r$ shows that they are quite degenerate. Recently  I.Vaisman
proposed the natural generalization of the geometric quantization
scheme in the case of Poisson manifolds. Though the r-matrix
structures we introduce in the context of the deformation
quantization  and they are quasiclassical limit of the quantum group
structures, it would be interesting to check the Vaisman condition in
this natural setting.

  Our goal is to verify the Vaisman condition of geometric
quantization for R-matrix type Poisson brackets  $\{\;,\;\}_{SD} +
\lambda \{\;,\;\}_{KKS}$.
  We are show that generalization of the quantization conditions is
failed in the case of hermitian symmetric spaces. We give two simple
examples when this condition is valid for the Poisson 2-matrix
structures.

In the first part of this work we introduce a Poisson pencil
connected with the standard modified {\it r}-matrix and formulate the
problem. In the second part we review  some basic results of the
theory of Poisson Lie groups and coset spaces. We also check the
conditions when the Poisson pencil is degenerate. In the third part
we discussed the Vaisman conditions for the geometric quantization of
this Poisson pencil.

\section{Sklyanin-Drinfeld brackets on Hermitian Symmetric Spaces.}

 {\bf Definition.}  {\sl A Lie group $G$ is {\em called a Poisson
Lie} group if it is a Poisson manifold such that the multiplication
$G\times G \rightarrow G$ is a morphism of Poisson manifolds, where
$G \times G$ is equipped with the product Poisson structure} .

  Every Poisson Lie structure on a semisimple connected or a compact
Lie group can be written in the coboundary form
\smallskip
\begin{equation}
\pi (g) = l_{g\ast}(r) - r_{g\ast}(r)
\end{equation}
\smallskip
 where ${\it r} \in \bigwedge ^{\otimes 2}{\g}$   satisfies the
modified Yang-Baxter equation, i.e.  $[[{\it r,r}]]$ is
$AdG$-invariant $\bf [2].$

{\bf Definition}. {\sl A Lie subgroup  $H$ of  a Poisson Lie group
$G$ which has its own Poisson Lie structure is called a {\em Poisson
Lie subgroup} if the inclusion  $i:H  \hookrightarrow  G$  is a
Poisson morphism.

 A coset space  $G/H$  whis a Poisson structure is called a {\em
Poisson  coset space} if the natural map  $G \rightarrow  G/H $ is a
Poisson one.}

 Let us assume that  $H$ is connected. There exists the Poisson
structure   on  $G/H $ if a subspace  $C^{\infty} (G/H)$  of
$H$-right invariant functions on $G$   is a subalgebra of Poisson
 algebra  $C^{\infty} (G) {\bf [9]}. $

 The intrinsic derivative of a Poisson tensor $\pi$
 \smallskip
 \begin{equation}
d_e : {\g}\rightarrow {\g}\wedge {\g}
\end{equation}
\smallskip
 given by  $X\rightarrow  (L_{\bar X}\pi ) (e)$  ,  where   $\bar X$
is any vector field on  $G$  with  $\bar X (e) = X$   define a Lie
algebra structure  $[ \;,\;]^{\ast} :\bigwedge ^{\otimes
2}{\g}^{\ast} \rightarrow  {\g}^{\ast}$ on  ${\g}^{\ast}$     dual
to  $d_e =  \delta$  and hence a Lie bialgebra structure on
${\g}\oplus {\g}^{\ast}$.

{\bf Proposition.  [9]} {\sl  $H $ is a Poisson Lie subgroup iff
${\h}^{\bot}$     is  an ideal   in  ${\g}^{\ast}$           ( where
${\h}$  is a Lie algebra
of H ). $G/H$ is a Poisson coset space  iff   ${\h}^{\bot}$   is a
subalgebra.}

   Let  ${\g}$   be a semisimple Lie algebra over  {\bf C} ,
${\g}_{0}$ be the same algebra considered  over {\bf R}.
  Let  ${\k}$ be a compact real form of  ${\g}$ , i.e.  ${\k} $  is a
fixed point  set of a standard  Chevalley  antiinvolution  $\sigma$
of   ${\g}$      :
\smallskip
\begin{equation}
\smallskip
 \sigma ( E_{\pm \alpha } )= -E_{\mp \alpha } , \sigma (H_{\alpha} )=
-H_{\alpha}, \sigma (\lambda X) = \bar {\lambda} \sigma (X)
\end{equation}

  where  ${ E_{\pm \alpha} , H_{\alpha} } $ is a Chevalley  basis of
${\g}$   with respect to a Cartan subalgebra   ${\h}$.

  Let  ${\h}_c$       be a Cartan subalgebra of  ${\k}$ and  ${\h}
={\bf C}{\h}_c$  . If we choose the system  $\Delta_{+}$  of positive
roots of  ${\h}$          in  ${\g}$      ,  then we have a
decomposition
\smallskip
\begin{equation}
{\g}_o = {\k}+{\ga}+{\n}_{+}
\end{equation}
\smallskip
 where  ${\n}_{+} = \sum_{\alpha \in \Delta_{+}} {\g}_{\alpha} $
,   ${\ga} = i{\h}_c$           is non-compact part of the Cartan
subalgebra    . This  additive decomposition  leads to the  {\em
Iwasawa  decomposition}  $K\times A \times N  \rightarrow  G_0.$

  Using the Iwasawa decomposition we  introduce  a  Lie bialgebra
structure  for the complex semisimple  Lie  algebra  ${\g}$      in
the following  way:
\smallskip
 \begin{equation}
\delta ( X ) = [ X\otimes 1 + 1\otimes X , {\it r} ] ,
\end{equation}
\smallskip
where  $\it r$  is the Drinfeld-Jimbo R-matrix.
\smallskip
\begin{equation}
 {\it r}  = \frac{i}{2} \sum_{\alpha \in \Delta _{+}}
E_{\alpha}\wedge E_{-\alpha}
\end{equation}
\smallskip
 In the compact case Lie bialgebra structure for ${\k}$ is a real
form of the standard Lie bialgebra structure for  ${\g}$      , given
by  r-matrix
\smallskip
\begin{equation}
{\it r}_o =  \frac{1}{4} \sum_{\alpha \in \Delta +}V_{\alpha} \wedge
W_{\alpha}
 \end{equation}
 \smallskip
  Let  $P$ be a parabolic subgroup in  $G$  . Every coadjoint orbit
of a compact group  $K$  for a fixed point  $x$          is naturally
isomorphic to a coset  $K/K_p$ ,  where $K_p = P\cap  K$.

  $K_p$ is a Poisson subgroup of the Poisson Lie group $K$. Hence  we
conclude  that  every  coadjoint orbit  $O$ in $K$ has the natural
Poisson structure  $\pi_{SD}$  given by the Poisson reduction.
The symplectic leaves of this structure are the $B$-orbits on $G/P
\simeq  K/K_p$   or Bruhat cells  $\bf [8]$.

 {\bf Definition}. {\sl We shall  call two Poisson structures {\em
compatible} if every linear combination of them is also a Poisson
structure. This means that   Schouten-Nijenhuis  bracket  of the
corresponding  bivector fields  is  equal  to  zero.}

  Let  $K$  be a  compact  real form  of  a simple  complex  Lie
group  $G$  and  $P$  be  a  parabolic  subgroup in  $G , K_p = P
\cap  K  , {\p} $ and ${\k}_p$  be Lie algebras  of  $P$  and  $K_p$.
Denote  by  $\Delta _p$    the  following  subset  of  positive
roots
\smallskip
\begin{equation}
\Delta _p = \{ \alpha \in \Delta _+ | E_{-\alpha} \not\in {\p} \}
\end{equation}
\smallskip
  It is usefull to remark  that  in the case  of  symmetric
decomposition
\smallskip
\begin{equation}
 {\k} = {\k}_x \oplus {\m}_x
\end{equation}
\smallskip
   is  only  one  simple  root $\alpha$      belongs  to $\Delta _p
{\bf [5]}.$

  Let  ${\cal O}_x \simeq K/K_p$  be  a  coadjoint   orbit  for  $k$
with  a  fixed  point  $x$ .

 {\bf Theorem.} {\sl $\bf [6]$   The  Kirillov-Kostant-Souriau $\{  ,
\}_{KKS}$   on  ${\cal O}_x$  is  compatible  with  the  Poisson  Lie
structure  $\{ , \}_{SD}$  induced  by  the  coboundary  Poisson  Lie
structure
on  $G$  if  and  only  if  this  orbit  is  a  hermitian  symmetric
space.  }.

\medskip

{\bf Proposition.} {\sl The structure $\pi_{\lambda}=\pi_{SD}+\lambda
\pi_{KKS}$ on a hermitian symmetric space is degenerate iff $\lambda
\in [-2,0]$.}

{\bf Proof. }

 Let us complete the ortonormal system
\smallskip
\begin{equation}
 \{ V_{\alpha},W_{\alpha} | \alpha \in \Delta_p\} \cap \{V_{\alpha}
,W_{\alpha} | \alpha \in \Delta_{+} \backslash \Delta_p \} \nonumber
\end{equation}
\smallskip
 to the ortonormal  basis with respect to the scalar product
$-\frac{1}{2} tr(\;,\;)$. The coadjoint action of $K$ preserves
Kartan-Killing form $tr(\;,\;)$ on ${\k}$ so we are abble to identify
\\
$Adg^{-1}=(Adg)^{\ast}=(Adg)^t$.

 Let us denote the corresponding matrixes of bivectors $r_o,r_p$ and
$r_{\lambda}$ , where
 \smallskip
 \begin{equation}
 r_{\lambda}=l_{g_{\ast}^{-1}}\pi_{\lambda}(g)=r_o-Adg^{-1}r_o+\lambda r_p
 \end{equation}
 \smallskip
 as $\bf R_o,R_p,R_{\lambda}$ respectively.

If $dim{\cal O}_x=2m$ then
\medskip
\begin{equation}
{\bf R_p}=\frac{1}{4}
\left( \begin{array}{ccc}
J_{2m} & 0 & \ldots \\
0 &\; 0 &\ldots \\
\vdots & \vdots & \ddots
\end{array} \right) ,
where
{\bf J_{2m}}=
\left( \begin{array}{ccc}
0 &1 & \ldots \\
-1 & 0  & \ldots \\
\vdots & \vdots & \ddots
\end{array} \right)
\end{equation}
\medskip
It is easy to show that the structure $\pi_{\lambda}(g)$ is
degenerated if and only if the left upper $2m \times 2m$-minor of
$\bf R_{\lambda}$ is degenerate. The rank of the left upper $2m
\times 2m$-minor is equal to the rank of $\bf16R_{\lambda}R_p$.\\
It's evident to see that
\smallskip
\begin{equation}
\bf det\left( {\rm 16}Adg^tR_oAdgR_p+(\lambda + {\rm 1}E_p\right)={\rm 0}
\end{equation}
\smallskip
where $\bf  E_p=-16R_oR_p=-16R_p^2$\\
\smallskip
We also see that
\smallskip
\begin{equation}
\bf det \left( {\rm 16}Adg^tR_oAdgR_p+(\lambda +{\rm1})E \right)={\rm
0}.
\end{equation}
\smallskip
This is a characteristic equation for matrix $\bf
{\rm16}Adg^tR_oAdR_p$.
\smallskip
Therefore
\smallskip
\begin{equation}
|\lambda +1| \leq || {\bf16Adg^tR_oAdgR_p }|| \leq 1 , where ||{\bf
A}||=supp_{||v||=1} \langle {\bf A}v , v \rangle.\\
\end{equation}
\smallskip
Now we have to prove the nesesary condition.

Let us consider a subgroup $K_{\alpha} \subset K$ with Lie subalgebra
${\k}_{\alpha}={\bf R}iH_{\alpha} \oplus {\bf R}V_{\alpha} \oplus
{\bf R}W_{\alpha}$, where $\alpha$ is a single root from $\Delta_p$.
It is well-known that such global subgroup exists.

The orbit of $K_{\alpha}$-action on ${\cal O}_x$ is naturally
isomorphic to the opdinary sphere. Moreover, symplectic leaves of
the induced structure on  a "small" orbit  coincide with the
intersection of symplectic leaves on ${\cal O}_p$ and this orbit ( if
we consider it as a Poisson coset space of the Poisson Lie subgroup
$K_{\alpha}$ . $\pi_{\lambda}$ on  $CP^1$ is degenerate for all
$\lambda \in [-2,0]$ . Thus we conclude that our structure
$\pi_{\lambda}$ is degenerate for all $\lambda \in [-2,0]$ too . $\bf
\Box$

\smallskip
Let $w$ be an element of the normalizator $N(H)$ of the maximal
length ih Wejl group $W(\Delta_+)$ . Then $Adw({\it r}_o)=-r_o$
because $Adw_{|h}:\Delta_{+} \rightarrow -\Delta_+$.

{\bf Proposition.} $\; \sl
l_{w^{\ast}}\pi_{\lambda}(g)=-\pi_{-(\lambda +2)}(g).$

{\bf Proof}
\smallskip
\begin{center}
$(l_{w^{\ast}})\pi_{\lambda}(g) = l_{g^{\ast}}(\pi_{\lambda}(w^{-1}g))
= l_{g^{\ast}}l^{-1}_{(w^{-1}g)^{\ast}}\pi_{\lambda}(w^{-1}g) =
l_{g^{\ast}}(r_o-Adg^{-1}(w^{-1}g)r_o+\lambda r_p) =
l_{g^{\ast}}(-r_o+Adg^{-1}r_o+(\lambda
+2)r_p)=-l_{g^{\ast}}(r_o-Adg^{-1}r_o-(\lambda
+2)r_p)=-\pi_{-(\lambda +2)}(g).
 \bf \Box$
\end{center}
\smallskip

That is why the topological properties ( i.e. topological structure
of symplectic leaves, quatnization condition ) of $\pi_{\lambda}$ and
$\pi_{-(\lambda +2)}$ are equivalent. If $\lambda =0$ then $\pi_o =
\pi_{DS}$ is equivalent to $\pi_{-\lambda}$.

\section{ On Geometric Quantization of R-type Poisson Pencil}

{\bf Definition.}{\sl The geometric quantization of a Poisson
structure is a map  of some subset of smooth functions on smooth
manifold $M$ to the set of hermitian operators which acts in the
space of the global cross sections of the hermitian line bundle over
$M$. Moreover, the following equation must be satisfied
\smallskip
\begin{equation}
\hat {\{f_1,f_2\}}=\frac{2\pi i}{\hbar }[\hat f_1,\hat f_2] ,
\end{equation}
\smallskip
where $f_1,f_2 \in C^{\infty}(M)$ and $\hat f$ is the image of $f$
(operator) , $\hbar$ is the Planck constant.}

In the symplectic case some fundamental results were obtained :
\smallskip
\begin{equation}
\hat f = f+\frac{\hbar}{2\pi i}\bigtriangledown_{c(f)}.
\end{equation}
\smallskip
Here  $c(f)$ is a vector field which acts as $c(f)g=\{f,g\}$ for all
$g\in C^{\infty}(M) , \bigtriangledown$ is a covariant derivative
which set a linear connection with the curvature form $\Omega$ and
$\Omega =\frac{1}{\hbar}\omega$ , where $\omega$ is a symplectic form
on $M$ . It was shown that $\Omega$ is a curvature form of some
complex line bundle if and only if the cohomology class of $\Omega$
is integer. The hermitian structure on the line bundle which endows
the space of global sections with Hilbert space structure exists and
is preserved iff this form $\Omega$ is real , i.e. $\omega\in\hbar
H^2(M,Z)$.

Some applications of this scheme were proposed by Vaisman $\bf [10]$
in a degenerated
case.

Let $M^m$ be a smooth manifold with a Poisson structure defined by a
bivector $\pi$.  Let $\cal A$ is a Poisson algebra, where ${\cal
A}=(C^{\infty}(M),\{\;,\;\})$.

We say that M permits a {\em prequantization} over a complex line
bundle $L\rightarrow M$ if the following identities are satisfied:
\smallskip
\begin{equation}
\hat{\{f_1,f_2\}}=\frac{2\pi i}{\hbar}[\hat f_1,\hat f_2] ,
\end{equation}
\smallskip
where $\hat f_j=f_j+\frac{\hbar}{2\pi i}\bigtriangledown_f$.

The differential operator $\bigtriangledown_f$ acts in the space of
global sections as
\smallskip
\begin{equation}
\bigtriangledown_fgs=\{f,g\}s+g\bigtriangledown_fs
\end{equation}
\smallskip
for all $f,g\in C^{\infty}(M) , s\in \Omega^o(L)$.

Let us consider an operator $\delta :\Upsilon^{\ast}(M)\rightarrow
\Upsilon^{\ast +1}(M)$ of Poisson differential in the following way:
\smallskip
\begin{equation}
\delta f=c(f) , \delta X=[[\pi ,X]] ,
\end{equation}
\smallskip
 where $[[\;,\;]]$ is the Schouten-Nijenhuis bracket, $f\in
C^{\infty}(M), X\in \Upsilon^{\ast}(M)$ is a polivector field.
It is easy to show that

${\cal P}(\lambda)(\alpha_1,...,\alpha_k)=(-1)^k\lambda ({\cal
P}(\alpha_1)^,...,{\cal P}(\alpha))$ for all
$\alpha_1,...,\alpha_k\in \Omega^1(M)$,
where $\beta ({\cal P}(\alpha))=\pi (\alpha ,\beta)\; \forall
  \alpha ,\beta \in \Omega^1(M) , \lambda \in \Omega^k(M)$ and ${\cal
P}(\lambda) \in \Upsilon^k(M)$ .\\
 $\cal P $ is a Hamiltonian operator defined by  $\pi$.\\
It was shown that the Jacobi condition for the Poisson bracket
$\{\;,\;\}$ is equivalent to the condition $\delta^2=o$ which is
satisfies iff $[[\pi ,\pi ]]=0$ . Therefore
$(\Upsilon^{\ast}(M),\delta )$ is a complex which is named the
Poisson complex with Poisson cohomologies {\bf [7].}

\medskip
\begin{equation}
H^k({\cal A})=\frac{Ker\{ \delta^k:\Upsilon^k(M)\rightarrow
\Upsilon^{k+1}(M)\}}{Im\{ \delta^{k-1}:\Upsilon^{k-1}(M)\rightarrow
\Upsilon^k(M)\}}
\end{equation}
\medskip

{\bf Theorem.} {\bf [10]}{\sl The Poisson manifold $(M,\pi )$ has a
quantization bundle if and only if there exist a vector field $X$ and
a closed 2-form  $\omega$ that represents an integer cohomology class
of $M$ , such that the relation
\smallskip
\begin{equation}
\pi +L_X\pi ={\cal P}(\omega) \label{main}
\end{equation}
\smallskip
holds.}

In that case a class $\frac{i}{2\pi}[c]\in H^2({\cal A})$ called the
{\em Poisson-Chern} class $c_1(\pi )$ is well-defined . It's image of
integral class from $H^2(M,Z)$ .

{\it Example  1.} Let $(p,q)\in R^2$ is standard symplectic space
with such brackets correlation:

$\{p,q\}=1$ .

We introduce a new structure

\smallskip
$\{f_1,f_2\}_{new}=\{f_1,h_1\}\{f_2,h_2\}-\{f_1,h_2\}\{f_2,h_1\}$ ,
where $h_1=p,h_2=pq$.

\smallskip
Here $\pi =p\frac{\partial}{\partial p}\wedge
\frac{\partial}{\partial q}$ and $\{p,q\}=p$. The Poisson structure
is degenerated if $p=0$.\\
Because of $H^{\ast }(R^2)=0$ we can get a trivial complex bundle
$R^2\times C \rightarrow R^2$ and the set of global cross sections
 is equal to $C^{\infty }(M)\otimes C$.
  That's why $\omega =0$.

Our equation $\pi +[[X.\pi ]]=0$ has a nontrivial solution
$X=q\frac{\partial }{\partial q}$. This bracket was also quatized by
means of the Drinfeld series $\bf [4].$

{\it Example 2.}  Let us consider an orbit of a nilpotent element in
${\g}^{\ast}$, where ${\g}={\it{sl}}(2,{\bf R})$.

  We shall identify  ${\g}^{\ast}$ with ${\g}$ using the scalar
product $tr<\;>$. Let us introduce
a coordinate system $(x_1 , x_2 ,x_3)$ in ${\g}$, such that every
matrix can be written in the form
\begin{equation}
\left(
\begin{array}{ccc}
x_2&x_1+x_3\\
x_1-x_3&-x_2
\end{array}
\right)
\end{equation}
 The equation $x_1^2+x_2^2=x_3^2$ corresponds to the singular orbit
of the nilpotent element.

In spherical coordinates
\smallskip
\begin{eqnarray*}
x_1&=&r\cos \phi \sin \theta \\
x_2&=&r\sin \phi \sin \theta \\
x_3&=&r\cos \theta
\end{eqnarray*}
\smallskip
this orbit is given by $\theta =\frac{\pi}{4}$.

 The Drinfeld-Jimbo R-matrix from $\bigwedge^2{\g}$ defines the {\it
r}-brackets
\smallskip
\begin{equation}
\{ r,\phi \}_{\it r}=-\frac{r}{2}\sin \phi
\end{equation}
\smallskip
which are degenerated if $\sin \phi =0$.

We construct a geometric prequantization in a trivial complex line
bundle as above:
\smallskip
\begin{equation}
\hat F =F+\frac{\hbar}{2\pi i} c(F) +X(F),
\end{equation}
\smallskip
where $X=-r\ln r\frac{\partial}{\partial r}$ is a solution of the
equation (\ref{main}).

 \bigskip

Let us study the geometric quantization of R-matrix type Poisson
family $\pi_{\lambda}$ on ${\cal O}_x$-orbit of a compact semisimple
Lie group $K$.

In particular case $K=SU(2) , {\cal O}_x=CP^1$ . $\pi_{\lambda}$ can
be written in the form
\smallskip
\begin{equation}
\pi_{\lambda}=-\frac{i}{2}(1+|z|^2)(\lambda +(\lambda
+2)|z|^2)\frac{\partial}{\partial z}\wedge \frac{\partial}{\partial
\bar z}.
\end{equation}
\smallskip
For $\lambda \in [-2,0]$ the Poisson bivector is degenerate on the
set
\smallskip
\begin{equation}
\Xi_{\lambda }=\{z| |z|^2=-\frac{\lambda}{\lambda +2}\} . \nonumber
\end{equation}
\smallskip
The geometric quantization condition is $\pi_{\lambda}+\delta X={\cal
P}\omega$, where $[\omega ]\in H^2(CP^1,Z)$

This is equal to ${\cal P}^{-1}(\pi_{\lambda})+d({\cal
P}^{-1}X)=\omega $ everywhere except $z\in \Xi_{\lambda}$ . In other
words
\smallskip
\begin{equation}
-\frac{2}{i}\frac{1}{1+|z|^2}\frac{1}{\lambda +(\lambda
+2)|z|^2}dz\wedge d\bar z = \omega -d\left (\frac{1}{\lambda
+(\lambda +2)|z|^2} \sigma \right ) ,
\end{equation}
\smallskip
 where $\sigma \in \Omega^1(CP^1).$

 We obtain using Stocks formula
\smallskip
\begin{equation}
\int \int_{|z|\geq \xi }{\cal P}^{-1}(\pi_{\lambda}) = \int
\int_{|z|\geq \xi }[-d\left ( \frac{1}{\lambda +(\lambda
+2)|z|^2}\right )\sigma +\omega]
\end{equation}
\smallskip
\begin{equation}
-2\pi \ln \frac{\lambda +(\lambda +2)\xi^2}{(\lambda
+2)(1+\xi^2)}=\phi (\xi )-\frac{1}{\lambda +(\lambda +2)\xi^2}\psi
(\xi ),
\end{equation}
\smallskip
where
\smallskip
\begin{equation}
\phi (\xi )=\int \int_{|z|\geq \xi}\omega \;\;\; and\;\;\;  \psi (\xi
)=\oint_{|z|=\xi } {\sigma} , \phi , \psi \in
C^{\infty}([\xi_o,+\infty ]) , \xi_o^2=-\frac{\lambda}{\lambda +2}
\end{equation}
\smallskip
If $\xi\rightarrow  +\xi_o$ the left part of the equation containes a
logarithmic singularity versus polar one of the right part . Hence
the geometric quantization condition {\em is not} satisfied contrary
to assumption. $\Box$

{\bf Proposition.}{\sl For all $\lambda \in [-2,0]$ the geometrical
condition for $\pi_{\lambda}$ is not  satisfied.}

{\bf Proof.} Let us consider again a subgroup $K_{\alpha}\subset K$
with Lie subalgebra ${\k}_{\alpha}=RiH_{\alpha}\oplus
RV_{\alpha}\oplus W_{\alpha} , \alpha \in\Delta_{+}$ and $\alpha$ is
single.

The orbit ${\cal O}_x$ of it's action is a Poisson submanifold with
respect to Poisson structure $\pi_{SD}$. Thus if $\lambda\not=0
\;\;\pi_{\lambda}(z,\bar z)$ is degenerate iff the corresponding
structure on "small" orbit with the same $\lambda$ is degenerate .
Hence the  equation $\pi_{\lambda}+\delta X={\cal P}(\omega) $ for
all nonsingular $z\in {\cal O}_{\alpha}$ can be reduced to the
equation $({\cal P}^{-1}\pi_{\lambda})+d( {\cal P}_{-1}(X))= \omega$
discussed above .

This equation on $CP^1$ has no solution . We show it using the
strsightforward calculations.

If  $\lambda =0$ then $\pi_{\lambda}=\pi_o=\pi_{SD}$ is degenerate
everywhere on the "small" orbit of $K_{\alpha}$ because it is a
Poisson submanifold. But  $l_{w^{\ast}} \pi_{\lambda}=-\pi_{-(\lambda
+2)}$ . So for $\lambda =0$ the structure {\em does not allow }
geometric quantization. $\Box $
\vskip 1.0cm
\section{Acknowledgments}

    I am deeply indepted to V.Rubtsov for pointing out
the problem and stimulating discussions. I would also to thank V.Fock
for useful suggestions and comments.

 This work was finished at the Institute of Theoretical Physics of
Uppsala University (Sweden).
I am thankful to professor A.Niemi for hospitality and excellent
conditions during my visit to Uppsala.

\end{document}